# Improving Profitability of a Color Production Line by Breaking Down Bottlenecks: A Computer Simulation Approach

Amir Jamali[1], Amirhossein Ranjbar[2], Ali Bozorgi-Amiri[3*]

## Abstract

Bottlenecks are one of the controversial issues in manufacturing companies. Managers and designers attempt to manage this trouble to improve efficiency in different ways. For example, expanding capacity is a prevalent solution to get rid of bottlenecks. In this paper, a color production line is chosen, which faces several challenges in its production line. This company attempts to distinguish and diminish the bottlenecks in the production line. The objective of this paper is to build a developed model of a production line to improve its profitability by breaking down its bottlenecks. Besides, the optimum number of operators with different utilizations is investigated in this paper. Furthermore, we considered the construction of new quality control in the station, which is the most time-wasting operation in the production line. The current study aims to apply computer simulation to examine the production line bottlenecks. In doing so, arena 14.00 software is used. Then the results are analyzed, and several managerial implications are presented.

[*] Corresponding author: Ali Bozorgi-Amiri; Email: alibozorgi@ut.ac.ir; Telefax: (+98) 21-8802-1067

[1, 2, 3,] School of Industrial Engineering, College of Engineering, University of Tehran, Tehran, Iran

## 1. Introduction and literature review

Supply chain management contains activities to manage materials and information flows between different echelons, such as retailers, manufacturing plants, and distribution centers. Supply chain problems cover a wide range of issues, from supply to sales, and analytical methods are usually used to tackle them. One of the most common methods is mathematical programming, which uses deterministic or stochastic models to resemble reality as much as possible. Still, for many of these problems it is impossible to employ analytic approaches. To address this issue, the simulation technique is employed when analytical methods are not useful to obtain performance evaluations. Although both of these methods seek optimization, there is a significant difference in their function. Analytical methods are highly dependant on the scenario that represents the empirical condition, and the answer is valid as long as the scenario does not change. Whereas, simulation examines a set of different scenarios to achieve the optimal solution. Therefore, simulation is a better approach to represent unstable conditions (Persson, 2011).

Recently, fabricating organizations are attempting to discover helpful solutions to diminish the common issues in the production line, for example, bottlenecks and holding up times. Organizations are endeavoring to get a competitive edge over their rivals by decreasing the additional expense and improving the production line. To reach these objectives, several solutions can be applied to manage various industrial issues that significantly impact the assembling system efficiency (Hutchison, 1991).

One of the critical applications of simulation is the design of production lines. Simulation helps managers compare the current state of the production line with the optimal condition in terms of performance. In the literature of simulation, using simulation in manufacturing systems is detailed in many papers (Kochar & Ma, 1989). Some simulation tecniques lead to more straightforward choices to manage the production line's issues (Terzi & Cavalieri, 2004). Adopting computer simulation in different domains such as production systems, emergency services, military industries, and construction projects has several advantages: financial benefits, enhancing resource utilization, reducing the cycle and flow time, and increasing throughput (Ivanov, 2018).

Computer simulation is a widely-used method to deal with production line issues. This method enables us to create theories and hypotheses for estimating future exercises dependent on the progressions in operational activities (Centeno, 1996). In fact, computer simulation has numerous applications in various fields, for example, improving the procedure, programming, booking, planning, production control, etc. One of the most valuable uses of simulation in the industry is improving efficiency by enhancing the asset usage, diminishing the process duration, and improving throughputs (Zahraee, Hatami, Yusof, Rohani, & Ziaei, 2013). ARENA software is a helpful simulation software widely applied by practitioners and researchers due to its capability to simulate uncertain scenarios (Memari, Zahraee, Anjomanshoae, & Rahim, 2013).

There is a growing body of research using simulation modeling due to its advantages over other approaches. The simulation method can be used to conduct a comparison between different alternatives and fault detection of existing systems. Persson (2011) presented a simulation tool based on Supply Chain Operation Reference (SCOR) model for simulating supply chain operations. The SCOR identifies best practices for supply chain operations and is a benchmark to improve supply chain performance. Yong et al. (2020) applied the simulation approach to predict the final capacities. They modeled the flow of air freights to evaluate the improvement suggestions. They contributed to the literature by examining the reused material flow in the recycling industry to balance capacity and demand. Kiani Mavi, Zarbakhshnia, & Khazraei (2018) modeled Tehran's Bus Rapid Transit (BRT) system with a discrete-event simulation. They investigated the effects of cost variables on the BRT performance using Arena 14 software to simulate Tehran's BRT network. To this end, they considered four scenarios and ranked them by multi-criteria decision-making methods. Line 1 of the BRT of Tehran is analyzed to show the applicability of their method.

Previous studies applied simulation methods to identify bottlenecks in their production line. They identified bottlenecks as the most crucial obstacle to product line optimization, and they used buffers and parallel machines to solve this problem. Baesler, Jahnsen, & Bio-bio (2003) formulated a discrete event simulation model to analyze the production system's performance by evaluating the utilization of each station. They simulated various alternatives to improve the process productivity. The authors also employed their model for a sawmill in Chile as a case study to show the validity of their model. Sharda & Akiya (2012) implemented simulation models for inventory management of perishable products. They asserted that simulation techniques could enhance inventory systems. Furthermore, they proved that their model outperforms previous models by simulating the inventory activities like scheduling, programming, and replenishment. Computer simulation is one of the useful methods to evaluate customer queues and waiting times. Shojaie, Haddadi, & Abdi (2012) proposed a model for assessing non-standard queues in compressed natural gas stations. The authors considered different scenarios. In each scenario, they simulated the stations' queues, and several insights were obtained for decreasing the length of queues and waiting time of customers.

## 2. Model Description

### 2.1. General description of model

This model simulates the production line of a color factory. The production line can be considered the most important part that affects all other supply chain activities. One way to design a production line is to use computer simulation and its capabilities to analyze its efficiency, which is used in this article. First, the current production line is simulated and analyzed. Several important reports such as waiting time queue, the number of materials in queue, and idle time and cost are provided. Then, the production line is detected imbalanced, and activities were planned to improve it. Generally, we seek to identify current product line faults and resolve them by making minor but effective changes and providing a better design of the production line.

### 2.2. Current Model

In this sub-section, we analyze the primary production line. Raw materials are first transported from the warehouse to the production line. In the first step, the resin is added to the caldron containing the raw material. Next, the caldron is transferred to a large mixer to make the paste color. The paste color is then transferred to one of the two existing permil machines to obtain the base color. The based color is brought to large mixers to add some solvents. By adding these solvents, the final product (processed color) is obtained. For quality control, several manufactured products are transferred to the laboratory to evaluate their quality. If the product does not have the required quality, it would be transferred back to the mixer to add the necessary materials. But if it meets the predefined standards, it will be weighed on the bascules to determine its weight and then enter the packing area for packaging. According to the customers' taste, there are three types of packaging at this station. One of them is all done by the machine and the other two by the operators. In the following, the production line of this manufacturing company is reenacted, and the production line's efficiency is improved. This organization is a leading producer of industrial and color industry. Since the items are produced considering the customer's request, the layout of the industrial facility depends on the job shop framework. The following figure illustrates the production line.

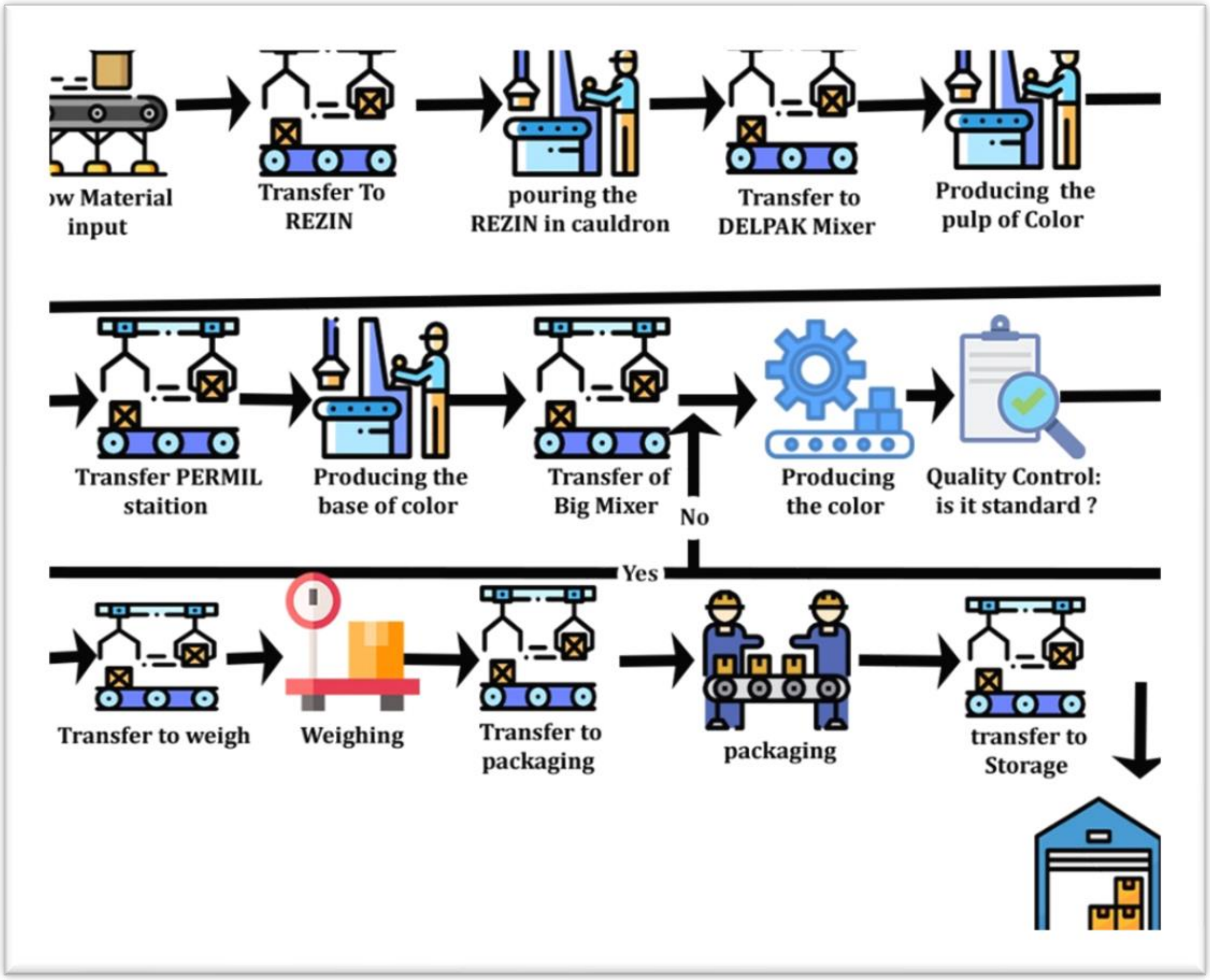

**Figure 1. The primary production line**

## 2.3. Problems

After simulating and examining the production line, imbalances in the production line were detected. The low utilization of a station indicates that it might be a bottleneck. The main effort of this article is to find bottlenecks and find a solution to reduce their harmful effects. The bottleneck station generates queues, increases the idle time of subsequent stations, and leads to other problems. Moreover, the low utilization rates of operators, the long queues of some machines, and the high idle time and cost of operators are recognized as the major problems of the primary model. In addition, color production has its own characteristics that prevent significant changes. For example, since some materials are chemical and quality control tests are so vulnerable, the laboratory distance to the production line must be maintained under all circumstances, which results in high shipping costs. The current model is simulated by the AENA 13.9 software for 50 replications. Each replication represents one shift (8 hours). The main results are shown in Table 1.

**Table 1. results of running current model**

| Items | Quantity | Unit |
|---|---|---|
| Average outputs | 38 | number |
| Average number in queue | 82.7 | number |
| Average waiting time in queue | 6.07 | hour |
| Average percentage of resource utilization | 60.3% | - |
| Value added cost | 7436 | dollar |
| Non-value added cost | 2035 | dollar |
| Busy cost | 9471 | dollar |
| Idle cost | 4418 | dollar |
| Total cost | 13889 | dollar |

### 2.4. Goals

The proposed model's main objective is to detect the bottlenecks and improve the efficiency of the production line. To achieve this goal, the current production line is simulated and analyzed, and the following developments were employed:

- determining the optimum number of operators for each station
- increasing the number of machines in bottlenecks
- implementing some operations simultaneously while transferring products in order to reduce the flow time
- checking the utilization rate of each machine and its operators
- investigating the probability of improving the current layout by changing the place of the quality control lab.

## 3. Developments

The utility rate of each machine is calculated, and the producing color station is recognized as the bottleneck of the production line. Therefore, a color production machine is added to the station with one operator. The analysis shows that utility rates of stations that are placed before the bottleneck are high, while those after the bottleneck are low. This indicates that the production line is imbalanced due to long queues and lack of buffers. By adding parallel machines, the efficiency rates of the machines and operators are balanced, and the idle cost is reduced.

Considering utility rates of operators, it is found that the unemployment rate of packaging operators is high. Hence, we determine the optimal number of operators in this station. The new operators are assigned to the added color production machine; hence, the number of production line operators has not changed. As a result of changes made, the packaging station's average utility rate rose from 20% to 80%.

By examining the product line, long transferring times are detected. To cover this issue, the mixing process is carried out by the mixing machines simultaneously with transportation to the next station. This resulted in saving time and energy while reducing the flow time.

Generally, the impacts of given developments can be seen by comparing pre and post-development tables. In addition to reducing costs and increasing productivity rates, the number of manufactured goods increased from 38 to 51 in the same production time. Furtheremore, improvements in the unit product's rate, the utility function of each operator, and the average waiting time in queues are depicted in figures 3, 4, and 5, respectively.

**Table 2. results of the developed model**

| **Items** | **quantity** | **unit** |
|---|---|---|
| Average outputs | 51 | number |
| Average number in queue | 73.15 | number |
| Average of waiting time in queue | 5.89 | hour |
| Average Percentage of Resource Utilization | 70.73% | - |
| Value added cost | 7547 | dollar |
| Non-value added cost | 1671 | dollar |
| Busy cost | 9218 | dollar |
| idle cost | 4663 | dollar |
| Total cost | 13881 | dollar |

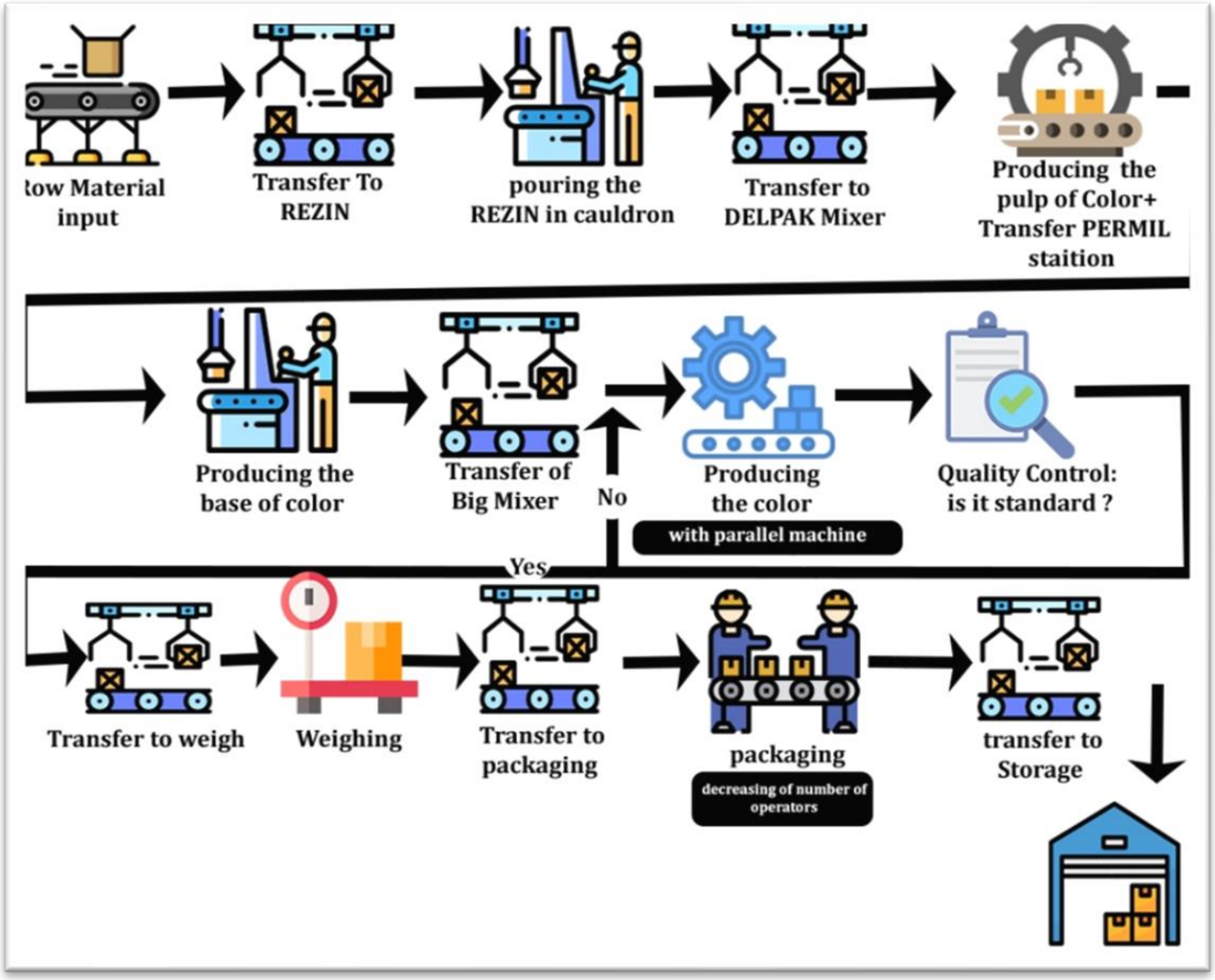

**Figure 2. The developed production line**

## 4. Result and Discussion

Figure 1 displays the current state map and the production line of the color factory. By simulating, we observed that the production line is not balanced, which causes lower efficiency and profit. To overcome this problem, first of all, we balance the company's production line. In order to do this, the remove bottleneck technique is performed. The bottleneck is identified as the most time-consuming task.

Figure 2 demonstrates the developed production line, which was drawn. In this figure, the map describes a developed model to balance the line. In this regard, some tasks are considered to be performed simultaneously. In addition, in some stations, parallel machines are considered to reduce cycle time. Furthermore, the model is endeavoring to find the optimum number of operators. In the following, the results of both models is compared. Figure 3 indicates that the number of output in the developed model increased by 38%. Also, the rate of production costs per product decreased by nearly 25%, which leads to more profit.

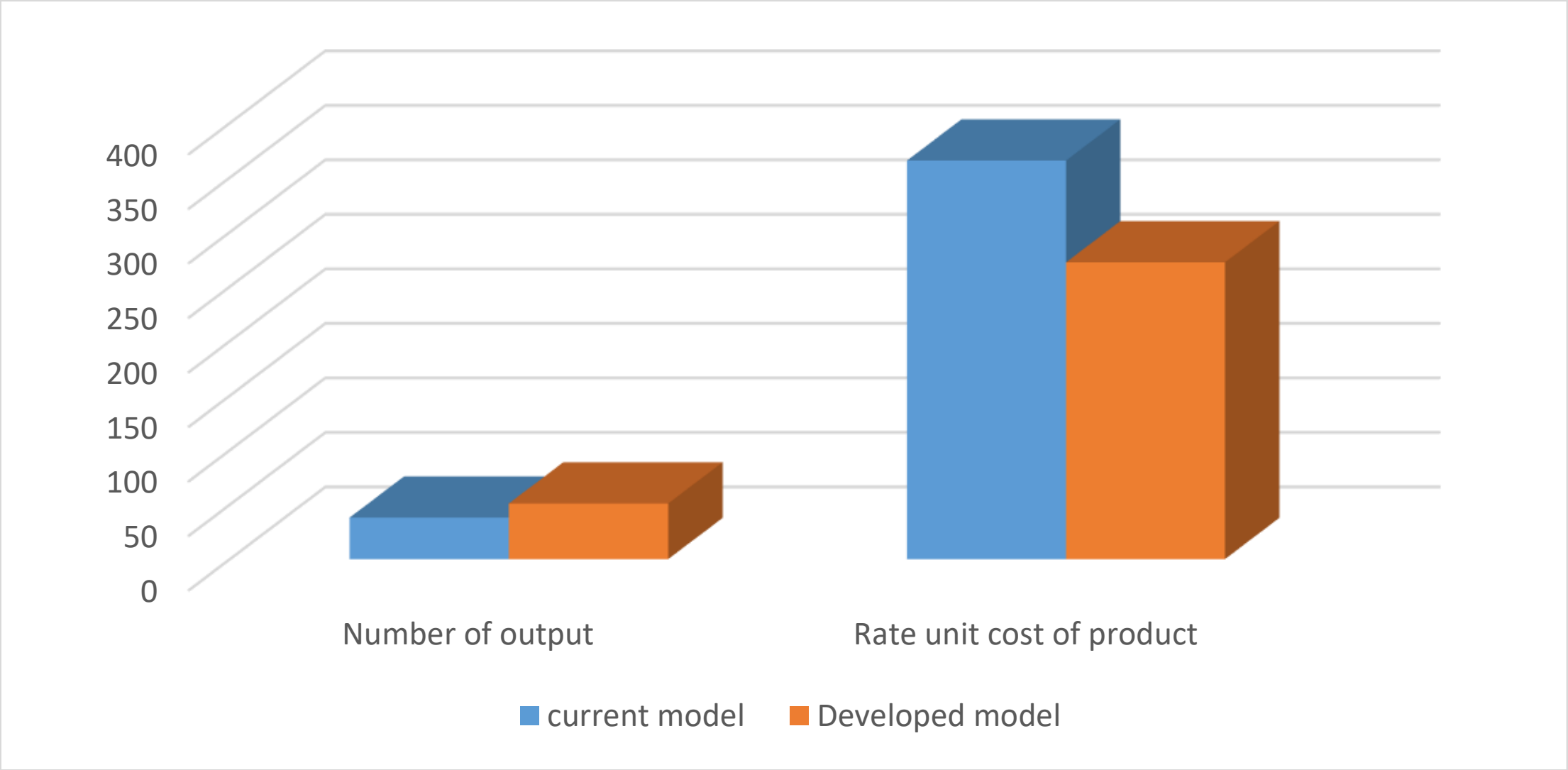

**Figure 3 .The rate of unit production in current and developed models**

Figure 4 shows that operators' performances are increased. Moreover, performance rates seem more equitable than the primary state because all operators work the same. Although it could make operators more tired and tedious, it should be noted the production rate and the amount of production are increased. Thus, managers can consider more rest time as an excellent strategy to increase employee satisfaction.

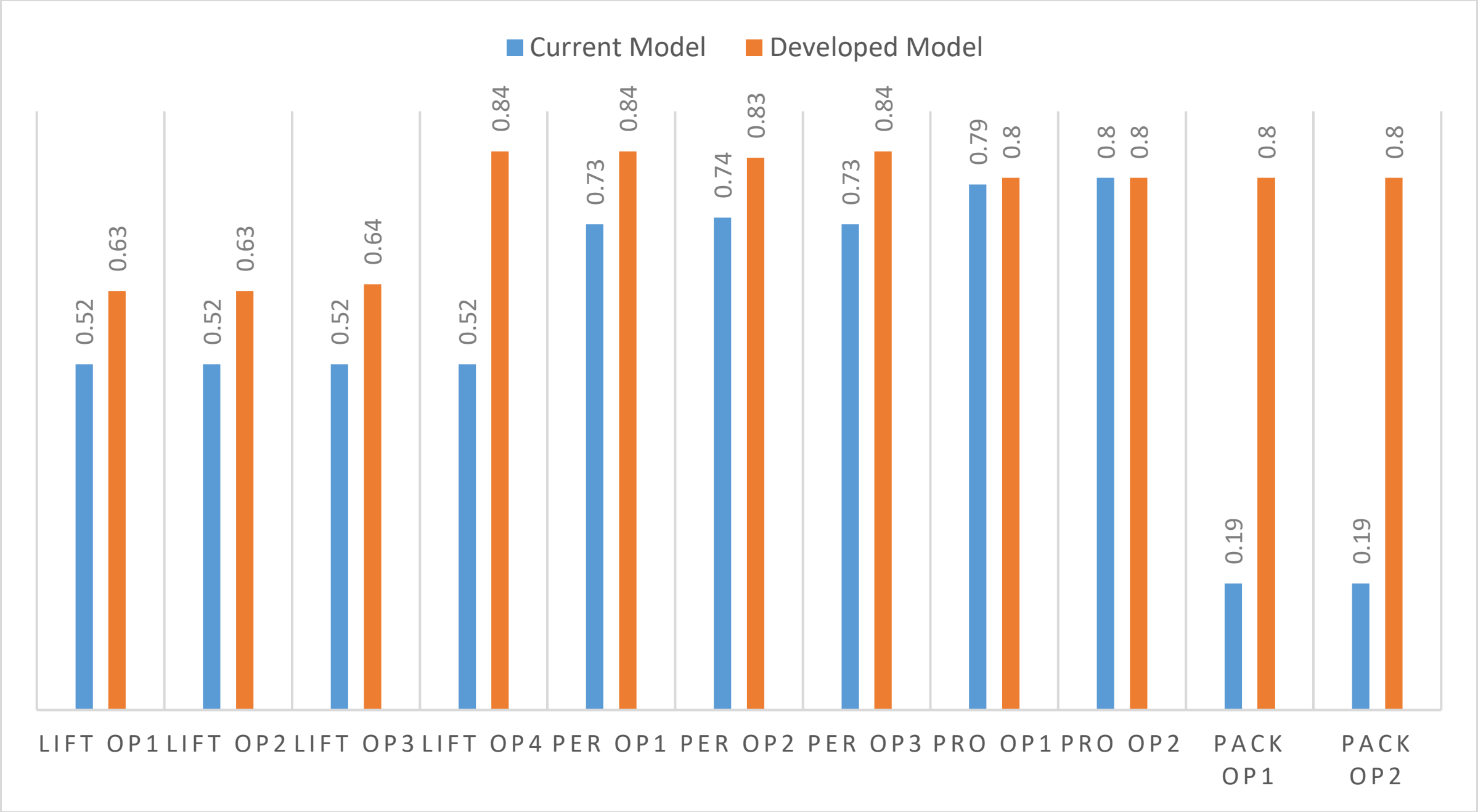

**Figure 4 The utility rate of each operator in the current and developed models**

Figure 5 shows that the production pulp of color and the color weighing stations have the most and least average time, respectively. Additionally, the difference between the maximum and minimum amount of average waiting time in the developed model is lower than the current model. Because in the developed model, the bottleneck waiting time is balanced.

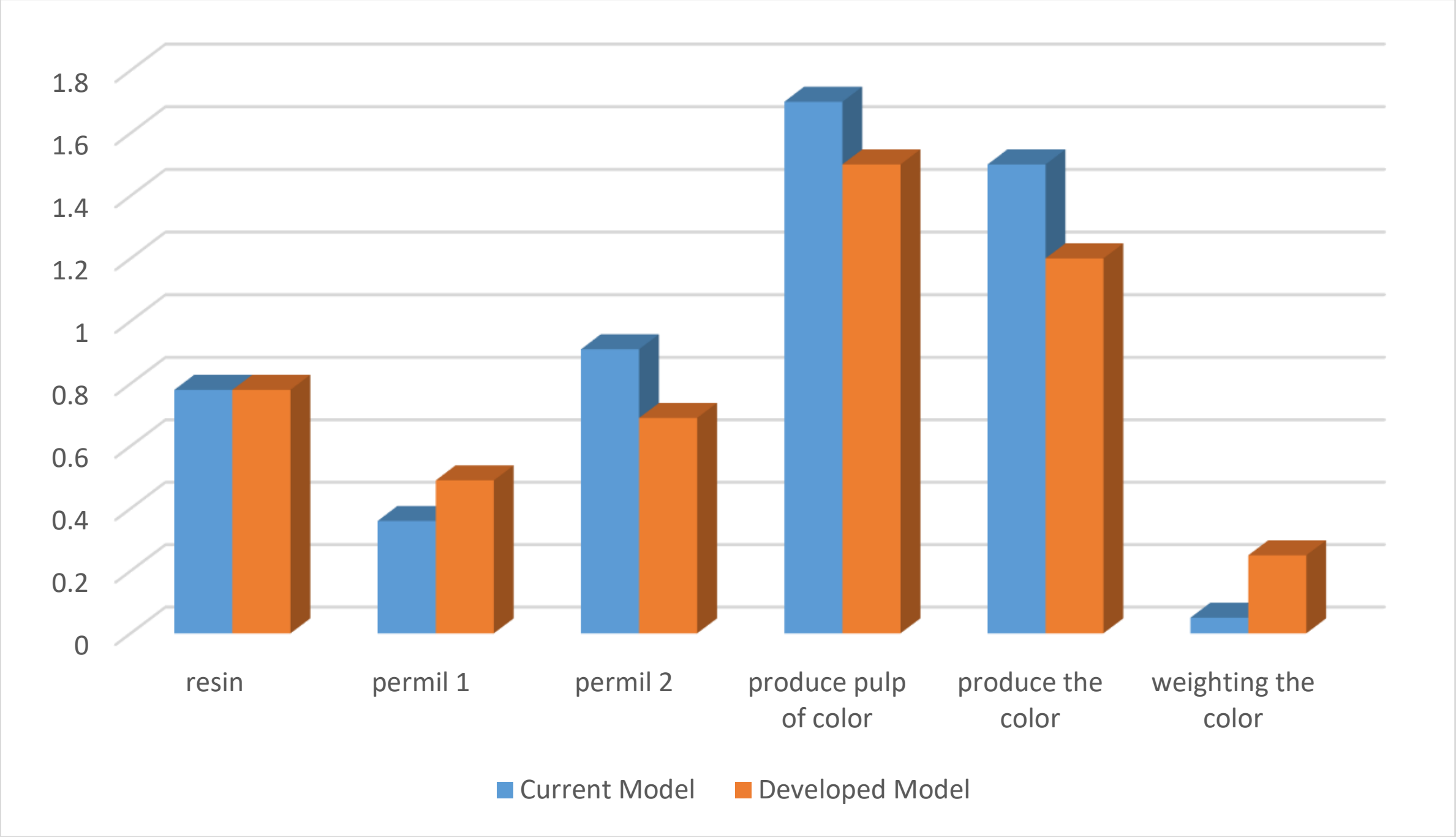

**Figure 5 Average waiting time in queues in the current and developed models**

## 5. conclusion

This study aims to present a better design of the production line by implementing minor changes in the color production line. Although specific features of color production make it difficult to improve it, we indicate how computer simulation can be employed to analyze improvement possibilities in a production line. Also, due to the long distance of the quality control lab to the production line, this part's cost and transferring time are very high. But the proposal to relocate the lab to the production line was rejected because of chemical materials and safety issues. Furthermore, a 10% increase in the productivity of the system is achieved by applying the discussed developments. Also, decreasing total cost and increasing produced goods are the other achievements of the above developments. In addition, relaxing some assumptions and restrictions of the study, considering failure events for machines, using maintenance and repairment schedulings can be interesting directions for future studies.

**Refrences**